\documentclass[12pt,a4paper]{amsart}
\usepackage{amscd,amssymb,amsopn,amsmath,amsthm,graphics,amsfonts,enumerate,verbatim,calc}

\headsep=1cm \numberwithin{equation}{section}
\def\bbR{\mathrm{I\!R}}

\def\dla{\mathcal{D}_{\hskip-2ptL}^*}

\def\lsq{\mathsf{[}}
\def\rsq{\mathsf{]}}
\def\hna{\hskip.2pt\widehat{\hskip-.2pt\nabla\hskip-1.6pt}\hskip1.6pt}
\def\dz{\mathcal{D}}
\def\hdz{\hskip.9pt\widehat{\hskip-.9pt\dz\hskip-.9pt}\hskip.9pt}
\def\hdp{\hskip.9pt\widehat{\hskip-.9pt\dz\hskip-.9pt}\hskip.9pt^\perp}
\def\hm{\hskip1.9pt\widehat{\hskip-1.9ptM\hskip-.2pt}\hskip.2pt}

\def\hg{\hat{g\hskip2pt}\hskip-1.3pt}

\def\w{^{\phantom i}}

\def\taw{{\tau\hskip-4.55pt\iota\hskip.6pt}}%
\def\hs{\hskip.7pt}
\def\hh{\hskip.4pt}
\def\hn{\hskip-.4pt}
\def\nh{\hskip-.7pt}
\def\nnh{\hskip-1pt}

\def\Gm{\Gamma}

\def\ve{\varepsilon}

\newtheorem{theorem}{Theorem}[section]
 
\newtheorem{lemma}[theorem]{Lemma}

\numberwithin{equation}{section}

\begin{document}
\pagenumbering{gobble}
\title{Compact Weyl-parallel manifolds}
\author{Andrzej Derdzinski}

\address{Department of Mathematics, The Ohio State University, 
Columbus, OH 43210, USA}
\email{andrzej@math.ohio-state.edu}

\begin{abstract}By 
ECS manifolds one means pseudo-Riemannian manifolds of dimensions $\,n\ge4\,$ 
which have parallel Weyl tensor, but not for one of the two obvious reasons:
conformal flatness or local symmetry.

As shown by Roter \cite{roter,derdzinski-roter-77}, they exist for every
$\,n\ge4$, and their metrics are always indefinite. The local
structure of ECS manifolds has been completely described
\cite{derdzinski-roter-09}.

Every ECS manifold has an invariant called rank, equal to 1 or 2. Known
examples of compact ECS manifolds
\cite{derdzinski-roter-10,derdzinski-terek-ne}, representing every dimension
$\,n\ge5$, are of rank 1. When $\,n\,$ is odd, some further, recently found
examples are locally homogeneous \cite{derdzinski-terek-cl}.

We outline the proof of the author's result, joint with Ivo Terek 
\cite{derdzinski-terek-tc}, which states that a compact rank-one ECS manifold,
if not locally homogeneous, replaced if necessary by a two-fold isometric
covering, must be the total space of a bundle over the circle.


\end{abstract}

\maketitle
\footnotesize{\noindent{Talk given at {\it The International Conference: 
Riemannian Geometry and Applications -- RIGA 2023}}\\
\phantom{Talk given at {\it The Intern}}Bucharest, September 22nd to 24th,
2023, online\\
\ \ \\
\ \ }

\section{The Ol\-szak distribution}\label{od}
\setcounter{equation}{0}
Given an ECS manifold $\,(M\nh,g)$, we define its {\it rank\/} 
to be the dimension $\,d\in\{1,2\}\,$ of its {\it 
Ol\-szak distribution\/} $\,\mathcal{D}$, which is a null parallel
distribution on $\,M\nh$. See \cite{olszak} and
\cite[p.\ 119]{derdzinski-roter-09}.

The sections of the Ol\-szak distribution $\,\mathcal{D}\,$ are the vector
fields $\,v\,$ such that 
$\,g(v,\,\cdot\,)\wedge[W(v'\nh,v''\nh,\,\cdot\,,\,\cdot\,)]=0$ for all
vector fields $\,v'\nh,v''\nnh$.

Lo\-rentz\-i\-an ECS manifolds have rank one: the Lo\-rentz signature limits
the dimensions of null distributions to at most $\,1$.

\section{Compact rank-one ECS manifolds}\label{cr}
\setcounter{equation}{0}
Examples of compact rank-one ECS manifolds have been found in all 
dimensions $\,n\ge5\,$ \cite{derdzinski-roter-10,derdzinski-terek-ne}. 

They are all geodesically complete, and none of them is locally homogeneous.
Recently \cite{derdzinski-terek-cl}, locally homogeneous examples (which are
necessarily incomplete) were exhibited in all {\it odd\/} dimensions
$\,n\ge5$.

It is an open question whether a compact ECS 
manifold may have rank 2, or be of dimension 4.

\section{The global structure theorems}\label{gs}
All known examples of compact ECS manifolds are dif\-feo\-mor\-phic to
nontrivial torus bundles over $\,S^1\nnh$, which reflects a general
principle -- here is our result \cite[Theorem A]{derdzinski-terek-ne}.
\begin{theorem}\label{maith}
Every non-locally-homogeneous compact rank-one ECS 
manifold, replaced if necessary by a two-fold isometric covering, is a
bundle over the circle, with the leaves of\/ $\,\mathcal{D}^\perp\nnh$
servng as the fibres.
\end{theorem}
The proof is outlined in Sections~\ref{so} --~\ref{sf}.

One needs the two-fold isometric covering in Theorem~\ref{maith} to make
$\,\mathcal{D}^\perp\nnh$ trans\-ver\-sal\-ly o\-ri\-ent\-able.

It is not known whether the conclusion of Theorem~\ref{maith} remains valid in
the locally homogeneous case, although it does then hold if 
$\,\mathcal{D}^\perp\nnh$ is also assumed to have at least one compact leaf.
The examples of locally homogeneous compact
rank-one ECS manifolds, constructed in \cite{derdzinski-terek-cl}, are bundles
over the circle.

A further result \cite[Theorem B]{derdzinski-terek-ne} pertains to universal
coverings of compact rank-one ECS manifolds:
\begin{theorem}\label{secnd}
The leaves of\/ $\,\hdp\hskip-2pt$ in 
the pseu\-do\hs-Riem\-ann\-i\-an universal covering space\/ 
$\,(\hm\nh,\hg)\,$ of 
any compact rank-one ECS manifold are the factor manifolds of a global product
decomposition of\/ $\,\hm\nh$.
\end{theorem}

Our notation uses hatted versions of symbols such as
$\,g,\dz,\,\mathcal{D}^\perp\nnh,\nabla\hs$ (the Le\-vi-Ci\-vi\-ta connection)
and $\,\mathrm{Ric}$, standing for objects in a given manifold $\,M\nh$, to
represent their analogs in the universal covering $\,\hm\nh$.

\section{The dichotomy property for foliations}\label{dp}
We refer to a co\-di\-men\-sion-one foliation $\,\mathcal{V}\hs$ on a 
manifold $\,M\,$ as having the {\it dichotomy property\/} when the following
condition is satisfied:

\vskip2pt
\noindent{\it Every compact leaf\/ $\,L\,$ of\/ $\,\mathcal{V}\hs$ has a 
neighborhood $\,\,U\,$ in $\,M\,$ such that\/ the leaves of $\,\mathcal{V}\hs$
intersecting\/} $\,\,U\nnh\smallsetminus\nh L$
\begin{enumerate}
\item[\ {\rm(i)}]{\it are either all noncompact, or}
\item[{\rm(ii)}]{\it they are all compact, and some neighborhood of\/
$\,L\,$ in\/ $\,M\,$
forming a union of compact leaves of $\,\mathcal{V}\hs$ may be
dif\-feo\-mor\-phic\-al\-ly identified with\/ $\,\bbR\nh\times L\,$
so as to make $\,\mathcal{V}\,$ the $\,L\,$ factor foliation.}
\end{enumerate} 

\section{Examples of the dichotomy property}\label{ed}
\setcounter{equation}{0}
Trans\-ver\-sal o\-ri\-enta\-bil\-i\-ty implies the dichotomy property when
\begin{enumerate}
\item[{\rm(a)}]$M\,$ and $\,\mathcal{V}\hs$ are real-an\-a\-lyt\-ic, or
\item[{\rm(b)}]$\mathcal{V}\hs$ has a finite number $\,r\ge0\,$ of compact
leaves.
\end{enumerate} 
For (a): any value of the leaf holonomy
representation, sending a neighborhood of $\,0\,$ in $\,\bbR\,$
real-an\-a\-lyt\-ic\-al\-ly into $\,\bbR$, must equal $\,\mathrm{Id}\,$ if it
agrees with $\,\mathrm{Id}\,$ on a nonconstant sequence tending to $\,0$.

Case (b) trivially follows be default. 
Examples of (b) include the Reeb foliation on $\,S^3\nh$,
while for any $\,r\ge0\,$ they obviously exist on $\,T^2$ and, consequently,
on $\,T^2\nh\times K$, with any compact manifold $\,K\nh$.

``Thickening'' a compact leaf $\,L\,$ satisfying (i) so as to
replace it with the closure of a prod\-uct-like
$\,\mathcal{V}\nh$-sat\-u\-rat\-ed neighborhood of $\,L$, one obtains a
foliation without the dichotomy property.

The dichotomy property easily follows in the case where $\,\mathcal{V}\hs$ is
the horizontal distribution of a flat linear connection in an orientable
(and hence trivial) real line bundle over a compact manifold $\,L$, with the
total space $\,M\nh$. 
Namely, the zero section $\,L\,$ is then a compact leaf, and depending on
whether the holonomy group of the connection is infinite or trivial, the
bundle has no global parallel sections except $\,L\,$ or, respectively, is
trivialized by them.

\section{The fibration lemma}\label{fl}
\setcounter{equation}{0}
\begin{lemma}\label{fibra}
Let a trans\-ver\-sal\-ly-o\-ri\-ent\-able 
co\-di\-men\-sion-one foliation\/ $\,\mathcal{V}\hs$ on a compact manifold\/ 
$\,M\,$ have the dichotomy property of Section\/~{\rm\ref{dp}}, and some
compact leaf\/ $\,L\,$ of\/ $\,\mathcal{V}\hs$ realize the option {\rm(ii)}, 
so that a prod\-uct-like\/ $\,\mathcal{V}\nh$-sat\-u\-rat\-ed neighborhood 
of\/ $\,L\,$ in\/ $\,M\,$ consists of compact leaves.

Then the leaves of\/  
$\,\mathcal{V}\hs$ are all compact, and 
constitute the fibres of a bundle projection\/ $\,M\nh\to S^1\nnh$.
\end{lemma}
This is \cite[Theorem\,4.1]{derdzinski-terek-tc}, and its proof uses the flow 
$\,\bbR\nh\times\nnh M\ni(\taw,x)\mapsto\phi\hs(\taw,x)\in M\,$ of a
$\,C^\infty\nnh$ vector field 
$\,\mathcal{V}\nh$. One fixes a point $\,z\,$ of a leaf satisfying condition
(ii) and applies a continuity argument to a maximal segment of the integral
curve $\,\taw\mapsto\phi\hs(\taw,z)\,$ which intersects compact leaves only.

Once we see that the maximal segment is defined on $\,(-\infty,\infty)$,
while any two leaves can obviously be joined by a piecewise $\,C^\infty$
curve, the smooth segments of which are integral curves of $\,C^\infty\nnh$
vector fields nowhere tangent to $\,\mathcal{V}\nh$, our claim becomes
reduced to 
a well-known exercise \cite[p.\,49]{moerdijk-mrcun}: a
trans\-ver\-sal\-ly-o\-ri\-ent\-able 
co\-di\-men\-sion-one foliation with compact leaves, on a compact manifold
$\,M\nh$, is tangent to the vertical distribution of a fibration
$\,M\to S^1\nnh$. 

\section{The local structure of rank-one ECS manifolds}\label{ls}
In coordinates $\,t,s,x\hh^i\nh$, where $\,i,j\in\{3,\dots,n\}$, the
following formula \cite{roter}, using constants $\,g_{ij}\w=g_{ji}\w$ and 
$\,a_{ij}\w=a_{ji}\w$, along with a function $\,f\,$ of the variable $\,t$, 
\begin{equation}\label{lsf}
\kappa\,dt^2\hn+\,dt\,ds\hs+\hs g_{ij}\w dx\hh^idx^j,\mathrm{\ with\
}\,\kappa=f\hn g_{ij}\w x\hh^ix^j\nh
+a_{ij}\w x\hh^ix^j
\end{equation}
defines a rank-one ECS metric if $\,f\,$ is nonconstant, 
$\,\det\hh[g_{ij}\w]\ne0=g\hh^{ij}a_{ij}\w$ and $\,[a_{ij}\w]\ne0$.

Conversely, at generic points (where $\,\mathrm{Ric}\,$ and
$\,\nabla\mathrm{Ric}\,$ are nonzero), any rank-one ECS metric has the above
form in suitable local coordinates. By lumping a rank-one ECS metrics together
with a special narrow class of 
locally symmetric ones, and allowing $\,f\,$ to possibly be constant, 
one gets rid of the genericity requirement
\cite[Theorem\,4.1]{derdzinski-roter-09}: (\ref{lsf}) always describes metrics
of this more general type, and all such have, in suitable local 
coordinates, the form (\ref{lsf}).

\section{Proof of the global structure theorem: four steps}\label{ps}
\setcounter{equation}{0}
\begin{enumerate}
\item[\ \ {\rm(I)}]We exhibit two functions $\,t,f:\hm\to\bbR\,$ 
\hbox{on the pseu\-do\hs-} Riem\-ann\-i\-an universal covering space 
$\,(\hm\nh,\hg)\,$ of a fixed compact rank-one ECS manifold
$\,(M\nh,g)\,$ in which $\,\mathcal{D}^\perp\nnh$ is 
trans\-ver\-sal\-ly o\-ri\-ent\-able, and introduce the 
space 
$\,\mathcal{S}\,$ of all continuous functions $\,\chi:\hm\to\bbR\,$ such that
the $\,1$-form $\,\chi\nnh\,dt\,$ is closed and projectable onto $\,M\nh$,
along with a linear operator 
$\,P:\mathcal{S}\to H^1\nh(M\nh,\hs\bbR)\,$ given by 
$\,P\chi=\lsq\chi\nnh\,dt\rsq$, where $\,\chi\nnh\,dt\,$ is treated as a
closed $\,1$-form on $\,M\nh$, and closedness of a continuous $\,1$-form
means its local exactness.
\item[\ {\rm(II)}]Using $\,t$, we prove the dichotomy
property of $\,\mathcal{D}^\perp\nnh$.
\item[{\rm(III)}]If $\,\dim\mathcal{S}<\infty$, local homogeneity follows.
\item[\ {\rm(IV)}]When $\,\dim\mathcal{S}=\infty$, the operator $\,P\,$ in (I)
is noninjective, and a nontrivial function in its kernel leads, via Sard's 
theorem, to a compact leaf $\,L\,$ of $\,\mathcal{D}^\perp\nnh$ realizing 
option (ii) of the dichotomy property, which allows us to use
Lemma~\ref{fibra}.
\end{enumerate}

\section{Step I: the functions $\,\,t\,\,$ and $\,\,f$}\label{so}
\setcounter{equation}{0}
We have $\,M\nh=\hm\nnh/\hh\Gm\hs$ 
for a group $\,\Gm\hs\,\approx\,\pi\nh_1\w\hn M\,$ 
acting on $\,\hm\,$ freely and properly dis\-con\-tin\-u\-ous\-ly via deck
transformations so as to preserve $\,\hg\,$ and the transversal orientation 
of $\,\hdp\nnh\nnh$.

The connection in $\,\hdz\,$ induced by the Le\-vi-Ci\-vi\-ta connection 
$\,\hna\,$ of $\,(\hm\nh,\hg)\,$ is flat: due to the lo\-cal-struc\-ture
formula (\ref{lsf}), $\,\hdz\,$ is spanned, locally, by the parallel gradient
$\,\hna\hn t$.

Simple connectivity of $\,\hm\,$ allows us to drop the word
`locally' and choose a global surjective function
$\,t\nnh:\nnh\hm\nh\to I\hs$
onto an open interval $\,I\nh\subseteq\bbR\,$ with parallel gradient
$\,\hna\hn t$, spanning $\,\hdz$.

This surjective function $\,t\nnh:\nnh\hm\nh\to I\hs$ 
is, clearly, unique up to af\-fine substitutions, and may be assumed, via
an af\-fine change, to coincide with the coordinate function $\,t\,$ 
in the lo\-cal-struc\-ture formula (\ref{lsf}).

Also, (\ref{lsf}) yields 
$\,\widehat{\hskip-1.7pt\mathrm{Ric}\hskip-1.7pt}
=(2-n)f\hs dt\otimes dt$, thus defining 
$\,f:\hm\to\bbR$, which is locally a function of $\,t$.

Consequently, any
$\,\gamma\in\Gm\hs$ gives rise to $\,q,p\in\bbR\,$ with $\,q>0$, such that,
for $\,(\hskip2.3pt)\hskip-2.2pt\dot{\phantom o}\nh=\hs d/dt$,
\begin{equation}\label{tge}
t\circ\gamma=qt+p,\ \ \gamma^*dt=q\,dt,\ \ f\circ\gamma=q^{-\nh2}f,\ \ 
\dot f\circ\gamma=q^{-3}f\nh.
\end{equation}
Closedness of a continuous $\,1$-form $\,\zeta$, such as
$\,\chi\nnh\,dt\in\mathcal{S}$,
means its being, locally, the differential of a $\,C^1$ function. The 
co\-ho\-mol\-o\-gy class 
$\,\lsq\hh\zeta\rsq\in H^1\nh(M\nh,\hs\bbR)
=\mathrm{Hom}\hh(\pi\nh_1\w\hn M\nh,\hs\bbR)\,$ then 
assigns to a homotopy class of piecewise $\,C^1$ loops at a fixed base point
the integral of $\,\zeta\,$ over a representative loop.

This results in a well-defined linear operator 
$\,P:\mathcal{S}\to H^1\nh(M\nh,\hs\bbR)$, where $\,P\chi
=\lsq\chi\nnh\,dt\rsq\,$ and $\,\chi\nnh\,dt\,$ is identified with the
projected $\,1$-form on $\,M\nh$.

\section{Step II: the dichotomy property of $\,\mathcal{D}^\perp$}\label{st}
\setcounter{equation}{0}
The normal bundle of a fixed compact leaf $\,L\,$ of $\,\mathcal{D}^\perp$ is 
canonically isomorphic, via $\,g$, to the line bundle $\,\dla$ over $\,L\,$
dual to $\,\mathcal{D}\hskip-2pt_L\w$ (the restriction of $\,\mathcal{D}\,$ to
$\,L$).

The horizontal distribution of the flat linear connection in $\,\dla$
arising from the one in the bundle $\,\mathcal{D}\,$ (spanned, locally, by the
parallel gradients $\,\nabla\hn t$) corresponds to the distribution
$\,\mathcal{D}^\perp$ under a suitable dif\-feo\-mor\-phic
identification $\,\varPsi\,$ of a neighborhood $\,\,U\hs$ of $\,L\,$ in
$\,M\,$ with a neighborhood $\,\,U'$ of the zero section $\,L\,$ in the line
bundle $\,\dla$.

The final paragraph of Section~\ref{ed}, slightly modified, then implies 
the dichotomy property.

We obtain the required dif\-feo\-mor\-phism $\,\,\varPsi\,$ using 
the flow $\,(\taw,x)\mapsto\phi\hs(\taw,x)\,$ of a fixed
smooth vector field on $\,M\nh$, nowhere tangent to $\,\mathcal{D}^\perp\nnh$,
and its lift $\,\hat\phi\,$ to $\,\hm\nh$. The resulting in\-te\-gral-curve
segments form the fibres of the tubular neighborhood $\,\,U\nh$, and along
these segments, pulled back to $\,\hm\nh$, denoting by $\,\pi:\hm\nh\to\ M\,$ 
the covering projection, we define $\,\varPsi\,$ by
\[
\varPsi(\phi\hs(\taw,x))\,
=\,[\hh t(\hat\phi\hs(\taw,y))-t(y)]\,\xi_y\w\nh\circ(d\pi\nh_y\w)^{-\nnh1}\,
\in\,\mathcal{D}_{\nnh x}^*\,\subseteq\,\dla\mathrm{,\ with\ }\,\pi(y)=x\hh,
\]
the parallel section $\,\xi\,$ of the line bundle $\,\hdz^*$ over 
$\,\hm\,$ being dual to $\,\hna\hn t\,$ in the sense that
$\,\xi(\hna\hn t)=1$. Hence $\,\varPsi\,$ sends local $\,t$-levels to local
sections parallel relative to the flat linear connection.

This construction is $\,\Gm\nh$-equi\-var\-i\-ant, and hence projects into
$\,M\nh$.

\section{Step III: the case $\,\,\dim\mathcal{S}<\infty$}\label{th}
\setcounter{equation}{0}
If $\,\dim\mathcal{S}\nh=m<\infty$, (\ref{tge}) and the final line 
of the text in Section~\ref{ls} give 
$\,|f|^{1/2},|\dot f|^{1/3}\nh\in\mathcal{S}$, while $\,\mathcal{S}\,$
is clearly closed under the $\,m$-ar\-gu\-ment operation 
$\,(\psi_1\w,\dots,\psi_m\w)\,\mapsto\,|\psi_1\w\dots\psi_m\w|^{1/m}$. 
Simple set-the\-o\-ret\-i\-cal reasons (see Appendix A) now
cause $\,|\dot f|^{1/3}$ to equal
a constant times multiple of $\,|f|^{1/2}\nh$, making $\,f$ globally a
function of $\,t$, of the form $\,f\nh=\ve\hh(t-b)^{-\nh2}$ with real
constants $\,\ve\ne0\,$ and $\,b$.

Combined with a result from algebraic geometry (Whitney's theorem), this
implies local homogeneity of $\,g$. See Appendix B.

\section{Step IV: the case $\,\,\dim\mathcal{S}=\infty$}\label{sf}
Now $\,P:\mathcal{S}\to H^1\nh(M\nh,\hs\bbR)\,$ is clearly noninjective.

Choosing
$\,\chi\in\mathcal{F}\smallsetminus\{0\}\,$ with $\,P\chi=0$, we see that 
$\,\chi\nnh\,dt$ projects onto an exact $\,1$-form on $\,M\nh$, that is, 
onto $\,d\mu\,$ for some (nonconstant) $\,C^1$ function $\,\mu:M\to\bbR$. 
As $\,\mathcal{D}^\perp\nnh=\hs\mathrm{Ker}\,dt\,$ on $\,\hm\nh$, this
$\,\mu\,$ is constant along $\,\mathcal{D}^\perp\nnh$.

Sard's theorem normally applies to
$\,C^k$ mappings from an $\,n$-man\-i\-fold
into an $\,m$-man\-i\-fold, for $\,k,n$ and $\,m\,$ with
$\,k\ge\mathrm{max}\hs(n-m+1,1)$,
guaranteeing that the critical values form a set of zero measure. 
In our case, $\,\mu:M\to\bbR$ is only of
class $\,C^1\nnh$, and $\,M\,$ can have any dimension $\,n\ge4$.

However, the conclusion of Sard's theorem remains valid here
\cite[Remark\,9.2]{derdzinski-terek-tc}, and so, due to
compactness of $\,M\nh$, the range $\mu(M)\,$ of $\,\mu\,$ 
contains an open interval formed by regular 
values of $\,\mu$.

In fact, $\,M\,$ is covered by
finitely many connected open sets $\,\,U\,$ each of which can be
dif\-feo\-mor\-phic\-al\-ly identified with an open set
$\,\,\widehat{U}\subseteq\hm\,$ such that the levels of 
$\,t:\widehat{U}\to\bbR\,$ are all connected. This turns $\,\mu$ restricted 
to $\,\,U\,$ into a function of $\,t$, allowing us to use Sard's
theorem as stated above for $\,k=n=m=1$.

Connected components $\,L\,$ of regular levels of $\,\mu$ clearly realize 
option (ii) of the dichotomy property.

\setcounter{section}{1}
\renewcommand{\thesection}{\Alph{section}}
\setcounter{theorem}{0}
\renewcommand{\thetheorem}{\thesection.\arabic{theorem}}
\section*{Appendix A}
\setcounter{equation}{0}
Here is an easy set-theoretical observation
\cite[Lemma\,3.3]{derdzinski-terek-tc}:
\begin{lemma}\label{setth}  
Let a vector space\/ $\,\mathcal{S}$ of functions\/
$\,X\to\bbR\,$ on a set $\,X\hs$ have a finite dimension 
$\,m>0$ and be closed both under the ab\-so\-lute-val\-ue operation\/ 
$\,\psi\mapsto|\psi|\,$ and under some\/ $\,m$-ar\-gu\-ment operation\/ 
$\,\varPi$ sending\/ $\,\psi_1\w,\dots,\psi_m\w$ 
to a function\/ $\,\varPi(\psi_1\w,\dots,\psi_m\w)\ge0\,$ having the same 
zeros as the product $\,\psi_1\w\nh\ldots\hs\psi_m\w$.

Then some basis of\/ $\,\mathcal{S}\,$ consists of nonnegative functions
with pairwise disjoint supports.
\end{lemma}
By `support' we mean {\it complement of the zero set.}

Applying the above lemma to our $\,\mathcal{S}\,$ we see that, on the set 
where $\,f\ne0$, the ratio $\,|\dot f|^{1/3}\nnh/|f|^{1/2}$ is locally
constant, which makes $\,|f|^{-\nnh1/2}$ (locally) linear as a function
of $\,t$.

Hence $\,f\ne0\,$ everywhere in $\,\hm\,$ (or else, at a boundary point
of the zero set of $\,f\nh$, the linear function $\,|f|^{-\nnh1/2}$ would
be unbounded on a bounded interval of the variable $\,t$).

Thus, $\,f\nh=\ve\hh(t-b)^{-\nh2}\nh$, as required, and $\,t\,$ has the range
$\,I\nh\subseteq\bbR\smallsetminus\{b\}$. Subjecting $\,t\,$ to an af\-fine
substitution, we may assume that $\,b=0\,$ and $\,I\nh\subseteq(0,\infty)$,
with $\,f\nh=\ve\hh t^{-\nh2}\nh$.

\setcounter{section}{2}
\renewcommand{\thesection}{\Alph{section}}
\setcounter{theorem}{0}
\renewcommand{\thetheorem}{\thesection.\arabic{theorem}}
\section*{Appendix B}
\setcounter{equation}{0}
Formula (\ref{lsf}) easily implies that the Le\-vi-Ci\-vi\-ta connection
$\,\hna\hs$ induces a flat connection in the quotient bundle 
$\,\hdp\hskip-2pt/\hdz\,$ over $\,\hm\nh$, with an
$\,(n-2)$-di\-men\-sion\-al pseu\-do\hs-Euclid\-e\-an space $\,V\nh$ of
parallel sections.

The (parallel) Weyl tensor naturally gives rise to a nonzero
trace\-less en\-do\-mor\-phism $\,A:V\to V$, represented by the matrix
$\,[\hs a_{ij}\w]\,$ in formula (\ref{lsf}).

Any $\,\gamma\in\Gm\nh$,
acting on $\,V\hs$ as a linear isometry $\,B$, pushes this $\,A\,$ forward
onto $\,B\nh AB^{-\nnh1}\nh=q^2\nnh A$, for $\,q$ related to $\,\gamma\,$ as
in (\ref{tge}). Due to compactness of $\,M\nh=\hm\nnh/\hh\Gm\nh$,
\begin{equation}\label{inf}
\mathrm{such\ }\,q\,\mathrm{\ arising\ from\ all\
}\,\gamma\in\Gm\hs\mathrm{\ form\
an\ infinite\ subset\ of\ }\,(0,\infty)\mathrm{,\ closed\ under
taking\ powers.}
\end{equation}
The set
$\,\mathcal{J}=\{(q,B)\in\bbR\times\mathrm{End}\,V:(BB^*\nh,B\nh AB^*) 
=(\mathrm{Id},q^2\nh\nnh A)\}\,$ is an algebraic variety 
in $\,\bbR\times\mathrm{End}\,V\nh$. By Whitney's classical result
\cite{whitney}, $\,\mathcal{J}\,$ has finitely many connected 
components, and 
hence so does the intersection $\,K=K'\nh\cap(0,\infty)$, for
the image $\,K'$ of $\,\mathcal{J}\,$ under the projection
$\,(q,B)\mapsto q$.

Thus, according to (\ref{inf}), $\,K=(0,\infty)$. 
Formula (\ref{lsf}), with $\,f\nh=\ve\hh t^{-\nh2}\nh$, now easily yields 
local homogeneity of $\,\hg$. Cf.\ also \cite{derdzinski}.




\end{document}